# Application of Extended Kalman Filter to Tactical Ballistic Missile Re-entry Problem

Subrata Bhowmik and Chandrani Sadhukhan

*Abstract*— **The objective is to investigate the advantages and performance of Extended Kalman Filter for the estimation of non-linear system where linearization takes place about a trajectory that was continually updated with the state estimates resulting from the measurement. Here tactile ballistic missile Re-entry problem is taken as a nonlinear system model and Extended Kalman Filter technique is used to estimate the positions and velocities at the X and Y direction at different values of ballistic coefficients. The result shows that the method gives better estimation with the increase of ballistic coefficient.**

## I. INTRODUCTION

THE effectiveness The effectiveness of the tactical ballistic missile (TBM) during the Persian Gulf War had a starling effect globally. Expertise developed during the War led to the availability of the anti tactical ballistic Missile (ATBM) for use in conventional air defence missile and, ultimately, to upgraded performance. The modern TBM was designed to maneuver in the terminal phase to escape tracking and interception. Therefore, on-line trajectory estimation of the TBM in the regular fire control sequence for the ATBM is highly desired in radar tracking, direct hit guidance, and early warning systems.

Normally trajectory estimation of a flying vehicle is investigated primarily in post analysis to identify the states and key parameters in available flight data measured using radar, satellites, and on board sensors. For a conventional ATBM air defence missile, engagement is limited in the re-entry phase since radar is the only measurement instrument used to sense the TBM. Thus, reconstructing the states and parameters of a re-entry vehicle is relatively difficult, particularly during maneuver.

Trajectory estimation of a manoeuvring TBM has received little attention. Chang et al. (1977) introduced a maneuvering re-entry vehicle filter with a newly defined augmented state vector. In the filter the augmented state vector contains position, velocity, and corresponding parameters for drag and manoeuvred forces. The Extended Kalman Filter has been proposed the task of estimation. The re-entry vehicle filter performance however is degraded if manoeuvring forces are absent. The related parameter estimation also relies heavily on the model and is assumed to be a constant Gaussian Markov process in the state equation [6]. Rabiner (1985) presented a Pontryagin filter based on the Pontryagin minimum principle, in which the unknown manoeuvring forces are treated as a control used to drive a vehicle's dynamics so that it will follow the noisy measured trajectory. Nevertheless, the Pontryagin filter is sensitive to parameter variation and initial conditions[3].

In TBM, manoeuvring is considered to be an extra lateral acceleration input, which rapidly changes a missile's position, velocity, and heading. Therefore, developing an input estimation approach to identify this lateral acceleration is highly desired. But in this case without developing input estimator only Extended Kalman Filter is applied on the system to estimate the state.

## II. EXTENDED KALMAN FILTER

In the Extended Kalman Filter linearization takes place about the filter's estimated trajectory. The partial derivatives are evaluated along a trajectory that has been updated with the filter's estimates. These in turn depend on the measurements, so the filter's gain sequence will depend on the sample measurement sequence realized on a particular run of the experiment. Thus the gain sequence is not predetermined by the process model assumptions as in the usual Kalman Filter [4].

A general analysis of the Extended Kalman Filter is difficult because of the feedback of the measurement sequence into the process model. The better trajectory is only better in a statistical sense. There is a chance that the updated trajectory will be poorer than the nominal one. In that event, the estimates may be poorer; this in turn leads the further error in the trajectory, which causes further errors in the estimates, and so forth, and so forth, leading to eventual divergence of the filter [5]. The net result is that the Extended Kalman Filter is a somewhat riskier filter than the regular linearized filter, especially in situations where the initial uncertainty and measurement errors are large. It may be better on the average than the linearized filter, but it is also more likely to diverge in unusual situations. However in Subrata Bhowmik is with the Department of Mechanical Engineering, National Institute of Technology, Rourkela, Orissa 769008 India (phone: 0661-2462524; e-mail: sbhowmik@nitrkl.ac.in).

Chandrani Sadhukhan (Roy) is with the Electrical Engineering Department, MCKV Institute of Technology, Howrah India (e-mail: future.iim@gmail.ac.in).

an Extended Kalman Filter it is usually more convenient to keep track of the total estimates rather than the incremental ones.

### III. DYNAMIC MODEL OF TBM IN RE-ENTRY PHASE

Considering a two-dimensional TBM in the re-entry phase and in a flat, non-rotating earth with constant gravity model as showed in fig 1. The distance that the TBM travel in this phase is shorter than the distance travel before re-entry. The setting for the TBM is a point mass and constant weight along a ballistic trajectory, where three types of significant forces act on the TBM. Among which gravity and drag are of two types. The advanced tactical ballistic missile usually steers at the right time in the terminal phase so as to escape radar tracking and missile interception. This steering may rapidly change the heading, velocity and acceleration of the TBM. The third force is introduced in a plane, which is perpendicular to the velocity vector if the TBM performs a maneuver. The force is induced by a maneuver and may be represented by components along $X_R$ and $Y_R$ which are relative to the acceleration terms in the equation of motion.

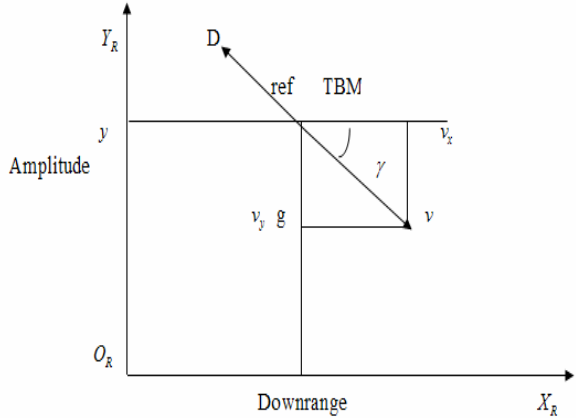

Fig. 1. Tactical ballistic missile geometry

The equations with a maneuver in radar co-ordinates $(O_R, X_R, Y_R)$ centered at a radar site can be represented as

$$\dot{v}_x = -\frac{Dg}{W}\cos\gamma + a_x = -\frac{\rho v^2}{2\beta}g\cos\gamma + a_x \quad (1)$$

$$\dot{v}_y = \frac{Dg}{W}\sin\gamma - g + a_y = \frac{\rho v^2}{2\beta}g\sin\gamma - g + a_y \quad (2)$$

$v_x$ and $v_y$ are the lateral velocities of the missile at X and Y direction, respectively. $\rho$ denotes the air density, $a_x$ and $a_y$ denote the lateral accelerations generated by a maneuver along $X_R$ and $Y_R$ direction respectively. The flight path angle $\gamma$ and ballistic coefficient $\beta$ are defined as

$$\gamma = \tan^{-1}(-\frac{v_y}{v_x})$$

$$\beta = \frac{W}{SC_{DO}}$$

Where $S$ and $C_{DO}$ represent the reference area and zero lift drag coefficient, respectively.
The air density is originally a function of altitude and should be considered into equation (1) and (2) because the altitude dramatically changes at flight at speeds over 100 km. the conventionally used approximation model for air density related to altitude is

$$\rho = .002378e^{-y/30000} \quad \text{for} \quad y < 30000\,ft$$
$$\rho = .0034e^{-y/22000} \quad \text{for} \quad y \geq 30000\,ft$$

Let the state vector be
$$X = [x_1\ x_2\ x_3\ x_4]^T = [x\ y\ v_x\ v_y]^T \quad (3)$$

The nonlinear state equation can be derived as
$$\dot{X} = F(X,t) + \phi u + I\xi \quad (4)$$

In equation (4), $\xi$ denote the process noise a vector with a variance of $Q$, $I$ stands for the identity matrix, $u$ is the lateral acceleration term of the TBM.

$$F(X,t) = \begin{bmatrix} x_3 \\ x_4 \\ -\frac{\rho}{2\beta}(x_3^2 + x_4^2)g\cos\gamma \\ \frac{\rho}{2\beta}(x_3^2 + x_4^2)g\sin\gamma - g \end{bmatrix}$$

$$\phi = \begin{bmatrix} 0 & 0 & 0 & 0 \\ 0 & 0 & 0 & 0 \\ 0 & 0 & 1 & 0 \\ 0 & 0 & 0 & 1 \end{bmatrix}$$

$$u = \begin{bmatrix} 0 & 0 & a_x & a_y \end{bmatrix}^T = \begin{bmatrix} 0 & 0 & u_3 & u_4 \end{bmatrix}^T$$

For a conventional air defence missile system, ground radar is the major instrument used to detect a TBM; it provides the position, velocity, and even the acceleration of tracked targets. The measurement equation is
$$Z = HX + \varepsilon \quad (5)$$

Where $\varepsilon$ represents the measurement noise vector and $H$ denotes the identity matrix. Equation (4) and (5) form the dynamic equations for the TBM with a maneuver after re-entry. Once all the states at a specific time instance are precisely known, the trajectory is reconstructed using a fixed-point smoother and n-step ahead predictor.

The Extended Kalman Filter is used as state estimator for nonlinear dynamic equations and is applied straight forwardly here in. Let $\varepsilon$ represents white and be normally distributed with zero mean and a variance of $R$. The predicted and updated state vectors from $t = n\Delta t$ to $t = (n+1)\Delta t$ under input $u_n$ at $t = n\Delta t$ are given by

$$\hat{X}_{n+1/n} = \phi_n \hat{X}_{n/n} + \varphi u_n \quad (6)$$

$$\hat{X}_{n+1/n+1} = \hat{X}_{n+1/n} + K_{n+1}(Z_{n+1} - H\hat{X}_{n+1/n}) \quad (7)$$

where $\Delta t$ is the sampling period, $Z_{n+1}$ denotes measurement at time $t = (n+1)\Delta t$, and the transition matrix

$$\phi_n = I + \frac{\partial F(X,t)}{\partial X}\bigg|_{X=\hat{X}_{n/n}} \Delta t$$

A linearization of $F(X,t)$ can be easily obtained

$$\frac{\partial F(X,t)}{\partial X}\bigg|_{X=\hat{X}_{n/n}} = \begin{bmatrix} 0 & 0 & 1 & 0 \\ 0 & 0 & 0 & 1 \\ 0 & f_{32} & f_{33} & f_{34} \\ 0 & f_{42} & f_{43} & f_{44} \end{bmatrix}_{X=\hat{X}_{n/n}}$$

$$f_{32} = \frac{\rho}{44000\beta} g(x_3^2 + x_4^2)\cos\gamma$$

$$f_{33} = -\frac{\rho}{2\beta} g(2x_3 \cos\gamma - x_4 \sin\gamma)$$

$$f_{34} = -\frac{\rho}{2\beta} g(2x_4 \cos\gamma + x_3 \sin\gamma)$$

$$f_{42} = -\frac{\rho}{44000\beta} g(x_3^2 + x_4^2)\sin\gamma$$

$$f_{43} = \frac{\rho}{2\beta} g(2x_3 \sin\gamma + x_4 \cos\gamma)$$

$$f_{44} = \frac{\rho}{2\beta} g(2x_4 \sin\gamma - x_3 \cos\gamma)$$

Kalman gain $K_{n+1}$ and Covariance matrix of $X_{n+1/n}$ and $\hat{X}_{n+1/n+1}$ i.e from $P_{n+1/n}$ and $P_{n+1/n+1}$ respectively [1,2,3] are

$$K_{n+1} = P_{n+1/n} H^T (H P_{n+1/n} H^T + R)^{-1}$$

$$P_{n+1/n} = \phi_n P_{n/n} \phi_n^t + \Gamma Q \Gamma$$

$$P_{n+1/n+1} = (I - K_{n+1} H) P_{n+1/n}$$

### IV. SIMULATION AND RESULTS

In this system the process noise variance and measurement noise variance are considered as .01 and .3 respectively and having zero mean. When Extended Kalman filtering method is applied on the system the simulation results are shown in the fig (2), fig (3), fig (4) and fig (5) respectively with taking the value of ballistic coefficient as $\beta=500$. The estimation of position X direction and Y direction, velocities X direction and Y direction are done. Again the model is tested with decreasing the value of ballistic coefficient as $\beta=300$. The results are shown in the fig (6), fig (7), and fig (8) and fig (9) respectively. The positional error for $\beta=500$ shows that takes time to estimate while taking $\beta=300$ the less time takes in this case. Same phenomena happen while estimating the position in Y direction. In case for velocity estimation technique in X and Y direction, it is found that with decrease of ballistic coefficient i.e. $\beta=300$ the time takes to estimate is less. For better performance analyze of the system four sets of plot is done by system is again run with increase the ballistic coefficient as $\beta=700$. The plots are shown in the fig (10), fig (11), fig (12) and fig (13). So the filter is performing well for this particular system.

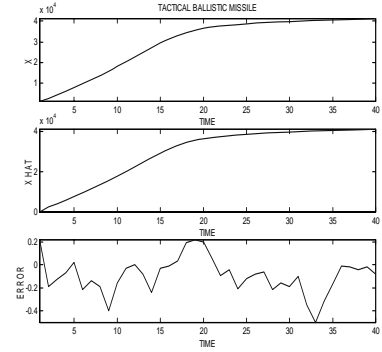

Fig. 2. Position estimation in X direction for $\beta=500$ Ib/ft$^2$

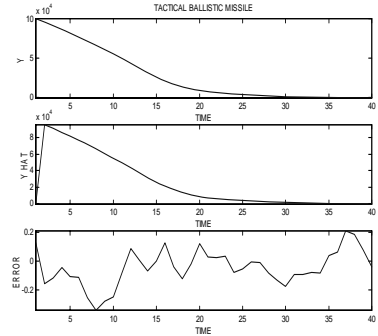

Fig. 3. Position estimation in Y direction for $\beta=500$ Ib/ft$^2$

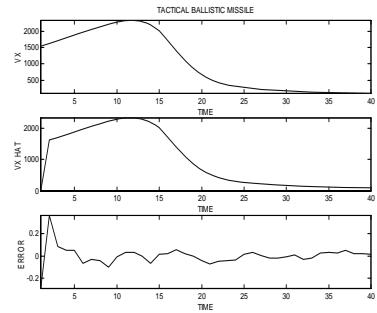

Fig. 4. Velocity estimation in X direction for $\beta=500$ Ib/ft$^2$

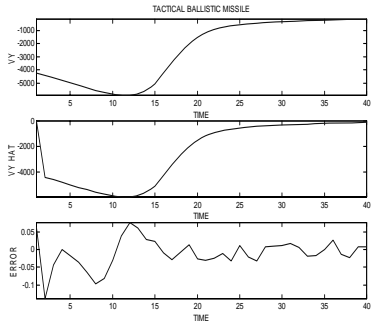

Fig. 5. Velocity estimation in Y direction for $\beta=500$ Ib/ft$^2$

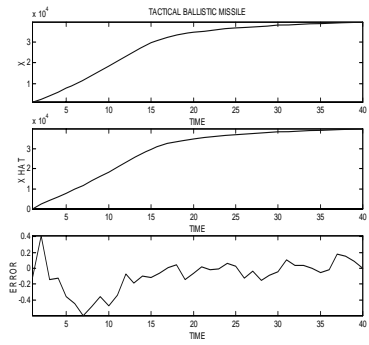

Fig. 9 Velocity estimation in Y direction for $\beta=700$ ib/ft$^2$

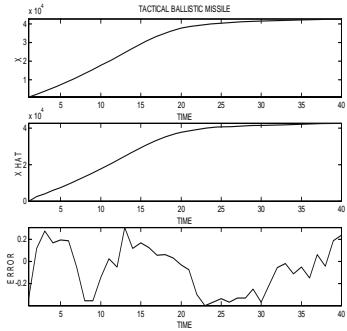

Fig. 6. Position estimation in X direction for $\beta=700$ ib/ft$^2$

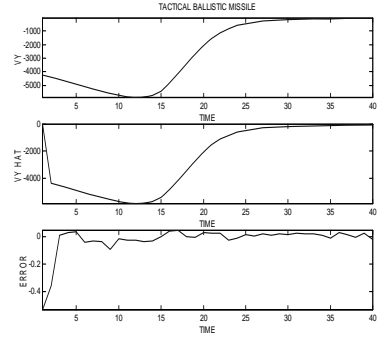

Fig.10 Position estimation in X direction for $\beta=300$ ib/ft$^2$

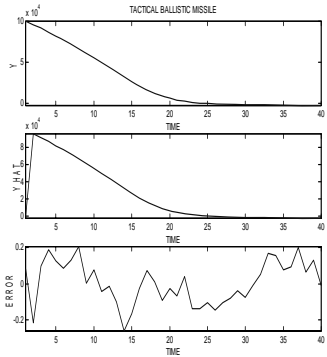

Fig. 7 Position estimation in Y direction for $\beta=700$ ib/ft$^2$

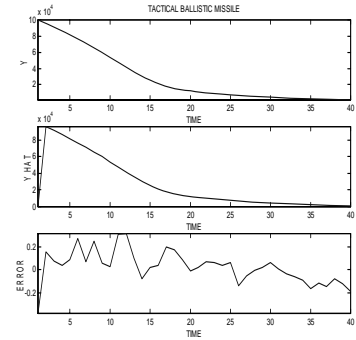

Fig.11 Position estimation in Y direction for $\beta=300$ ib/ft$^2$

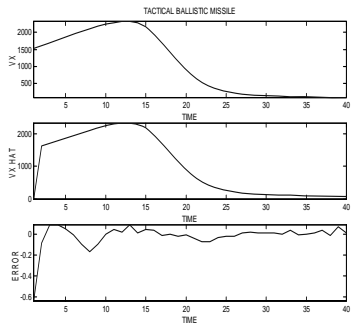

Fig. 8 Velocity estimation in X direction for $\beta=700$ ib/ft$^2$

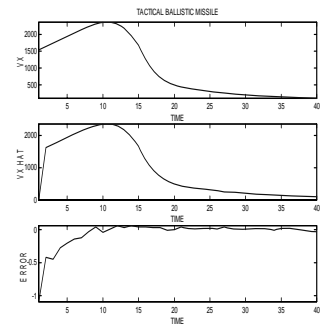

Fig.12 Velocity estimation in X direction for $\beta=300$ ib/ft$^2$

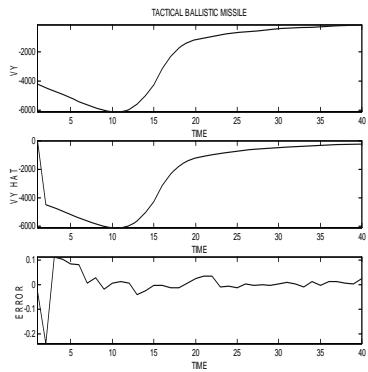

Fig.13 Velocity estimation in Y direction for $\beta=300$ ib/ft$^2$

V. CONCLISION

In this Paper Application of Extended Kalman Filter to the estimation of Nonlinear system is discussed. Tactile Ballistic Missile Re-entry problem is a very important problem in defence research and it is selected as a test model. From the results it is clearly shown that better performance is done with the increase of ballistic coefficient. The error is reduced with the increase of ballistic coefficient.

REFERENCES

[1] S Bhowmik, C Roy, "*Comparison of Estimation Techniques using Kalman Filter and Grid based Filter*", Accepted for publication in Proceeding of International Conference on advances in Control and Design optimisation in Dynamic systems, IISc,Bangalore,01-02 Feb,2007
[2] S Bhowmik, C Roy, "*Comparison of Estimation Techniques using Kalman Filter and Grid based Filter for Linear and Non-linear systems* ", Accepted for publication in Proceeding of  International Conference on Computing : Theory and applications , ISI, Kolkata , 04-05March ,2007
[3] C Roy , *Comparision of Kalman Filter and Grid based filter for Non Linear systems* , M.Tech Thesis , Jadavpur University,  Kolkata , 2003
[4] H Kushner, Report on   "*Approximation Methods for Nonlinear Filtering*" , Brown University 1993
[5] W   Kliemann, "*Nonlinear Filtering approaches to Multi target Tracking* " Journal of Signal Processing,Vol12,pp-123-128, July, 1997
[6] L.R Rabiner, "*A Tutorial on Hidden Markov Models and selected Application in speech Recognition*". Proc of the IEEE Conference on Signal Processing, Sydney, 1989